\documentclass[a4paper,11pt]{article}
\usepackage[latin1]{inputenc}
\usepackage{amsfonts}
\usepackage{amsmath}
\usepackage{amssymb}
\usepackage{amsthm}
\usepackage[T1]{fontenc}
\usepackage[dvips]{graphicx}
\usepackage{color}
\usepackage{epsfig}
\theoremstyle{plain}
\newtheorem{Thm}{Theorem}
\newtheorem{Lem}{Lemma}
\newtheorem{Def}{Definition}

\newtheorem*{Rem*}{\textsc{Remark}}
\newtheorem*{Lem*}{\textsc{Lemma}}
\newtheorem*{Cor*}{\textsc{Corollary}}
\newtheorem*{Con*}{\textsc{Conjecture}}
\newcommand{\e}{\varepsilon}
\newcommand{\bb}[1]{\mathbb{ #1 }}
\newcommand{\av}[1]{\left| #1 \right|}

\title{Poles of Int\'egrale Tritronqu\'ee and Anharmonic Oscillators. Asymptotic localization from
WKB analysis}
\author{Davide Masoero \thanks{E-mail address: masoero@sissa.it; Tel: +39 0403787501; Fax: +39 0403787528}\\ Mathematical Physics Sector, SISSA, via Beirut 2, 34151, Trieste}
\date{}
\begin{document}

\maketitle

\abstract{Poles of integrale tritronquee are in bijection with cubic oscillators
that admit the simultaneous solutions of two quantization conditions. We show that the poles
are well approximated by solutions of a pair of Bohr-Sommerfeld quantization conditions
(the Bohr-Sommerfeld-Boutroux system): the distance between a pole and the corresponding solution
of the Bohr-Sommerfeld-Boutroux system vanishes asymptotically.}

\section{Statement of Main Result}

In a previous paper \cite{piwkb}, the author studied the distribution of poles
of solutions of the the first Painlev\'e equation
\begin{equation*}
y''= 6y^2 -z \, , \; z \in \mathbb{C} \quad ,
\end{equation*}
with a particular attention to the poles of the int\'egrale tritronqu\'ee. This is the
unique solution of P-I with the following asymptotic behaviour at infinity

\begin{equation*}
y(z) \sim - \sqrt{\frac{z}{6}}, \quad \mbox{if} \quad |\arg z| <\frac{4 \pi}{5} \; .
\end{equation*}

The problem of computing the poles of the tritronqu\'ee solution was mapped to a pair of spectral problems for
the cubic anharmonic oscillator. More precisely, it was shown that a point $a \in \mathbb{C}$
is a pole of the tritronqu\'ee solution if and only exists $b \in \mathbb{C}$ such that the following
Schr\"odinger equation 
\begin{equation}\label{eqn:schroedinger}
\frac{d^2\psi(\lambda)}{d\lambda^2}= V(\lambda;a,b) \psi(\lambda)\; , \quad
V(\lambda;a,b)=4 \lambda^3 -2 a \lambda -28 b \;.
\end{equation}
admits the simultaneous solutions of two different quantization conditions.

Using a suitable complex WKB method, the author studied this pair of quantization conditions.
He derived a system of two equations, the Bohr-Sommerfeld-Boutroux (B-S-B) system,
whose solutions describe approximately the distribution of the poles.

\begin{figure}[htbp]
\begin{center}
\input{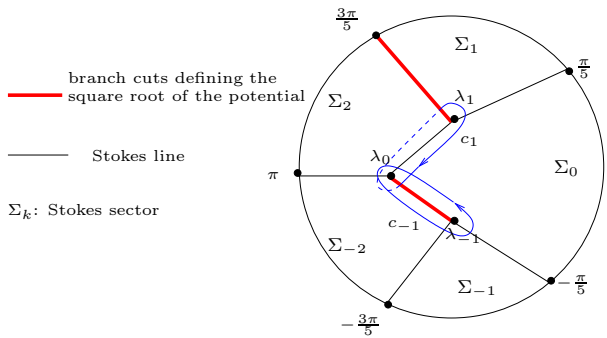}
 
\end{center}

\caption{Riemann surface $\mu^2=V(\lambda;a,b)$}
\label{figure:riemann}
\end{figure}

We say that $(a,b)$ satisfies the Bohr-Sommerfeld-Boutroux (for the precise
definition see Definition \ref{def:bsb} below) system if
\begin{eqnarray}\nonumber
 \oint_{c_{- 1}}\sqrt{V(\lambda;a,b)}d\lambda &=& i\pi (2n-1) \; , \\
\label{sys:bsb} \\ \nonumber
\oint_{c_{1}}\sqrt{V(\lambda;a,b)}d\lambda &=& i\pi (2m-1)\; .
\end{eqnarray}
Here $m,n \in \bb{N}-0$ are called quantum numbers and the cycles $c_{\pm1}$
are depicted in Figure \ref{figure:riemann}.

For any pair of quantum numbers there is one and only one solution to the
Bohr-Sommerfeld-Boutroux system; this is proven for example in \cite{kapaev03}.

Solutions of B-S-B system have naturally a multiplicative structure.

\begin{Def}\label{def:primitive}
Let $(a^*, b^*)$ be a solution of the B-S-B system with quantum numbers $n,m$ such that $2n-1$ and $2m -1$
are coprime. We call $(a^*,b^*)$ a {\rm primitive} solution of the system and denote it $(a^q,b^q)$, where $q=\frac{2n-1}{2m-1} \in \bb{Q}$. Due to Lemma \ref{lem:chiderivatives} below, we have that $$(a_k^q,b_k^q)=((2k
+1)^{\frac{4}{5}}a^q,(2k+1)^{\frac{6}{5}}b^q) , k \in \bb{N} \; ,$$ is another solution of the B-S-B system.
We call it a {\rm descendant} solution. We call $\left\lbrace (a_k^q,b_k^q)
\right\rbrace_{k \in \bb{N}}$ the q-sequence of solutions.
\end{Def}
In \cite{piwkb} it is shown that the sequence of real solutions of B-S-B system is the 1-sequence of solutions.
The real primitive solution is computed numerically as $a^1 \cong -2,34, b^1 \cong -0,064 $.

In the present paper we prove that any q-sequence
approximates a sequence of poles of the tritronqu\'ee solution.
The error between the pole and its WKB estimate is of order $(2k+1)^{-\frac{6}{5}}$ (see Theorem
\ref{thm:aapprox} below).

\begin{Def}\label{def:disc}
We denote $D_{\epsilon}(a')=\left\lbrace\av{a-a'} <\e , \, \e \neq 0\right\rbrace$.
\end{Def}

The main results of the present paper is the following

\begin{Thm}[Main Theorem]\label{thm:aapprox}
Let $\e$ be an arbitrary positive number. If $\frac{1}{5} <\alpha < \frac{6}{5}$,
then it exists $K \in \bb{N}^*$ such that for any $k\geq K$
inside the disc $D_{k^{-\alpha}\e}(a^q_k)$ there is one and only one pole
of the int\'egrale tritronqu\'ee.
\end{Thm}

The rest of the paper is devoted to the proof of the theorem.

\section{Proof}

\subsection{Multidimensional Rouche Theorem}

The main technical tool of the proof is the following generalization of the classical Rouch\'e theorem.

\begin{Thm}[\cite{aizenberg}]
Let $D,E$ be bounded domains in $\bb{C}^n $, $ \overline{D} \subset E$ , and let $f(z) , g(z)$ be holomorphic maps $ E \to \bb{C}^n$ such that
\begin{itemize}
 \item $f(z) \neq 0 ,\; \forall z \in \partial D$,
 \item $\av{g(z)} < \av{f(z)}, \;\forall z \in \partial D$, 
\end{itemize}
then $w(z)=f(z)+g(z)$ and $f(z)$ have the same number (counted with multiplicities)
of zeroes inside $D$. Here $\av{f(z)}$ is any
norm on $\bb{C}^n$.
\end{Thm}

\subsection{Monodromy of Schr\"odinger Equation}

Poles of int\'egrale tritronque\'ee are in bijections with the simultaneous solutions of two eigenvalues
problems for the cubic anharmonic oscillator \cite{piwkb}.
Below we recall the basics of anharmonic oscillators theory;
all the details can be found in \cite{piwkb}.

Fix $k \in \mathbb{Z}_5 = \left\lbrace -2, \dots,2 \right\rbrace$
and the branch of $\lambda^{\frac{1}{2}}$ in such a way that
${\rm Re}{\lambda^{\frac{5}{2}}} \to + \infty$ as
$|\lambda| \to \infty,\arg{\lambda}=\frac{2\pi k}{5}$. Then there exists
a unique solution $\psi_k(\lambda)$ of equation (\ref{eqn:schroedinger}) such that
\begin{equation}\label{def:psik}
\lim_{\lambda \to \infty, \av{\lambda -\frac{2 \pi k}{5}} < \frac{3 \pi }{5} -\varepsilon}
\lambda^{\frac{3}{4}}e^{+\frac{4}{5}\lambda^{\frac{5}{2}} - \frac{1}{2}
a\lambda^{\frac{1}{2}}}\psi_k(\lambda;a,b)=1 \; .
\end{equation}

For any pair of functions $\psi_l, \psi_{l+2}$, we  call \begin{equation}\label{def:wk}
 w_k(l,l+2)=\lim_{\substack{\lambda\ \to \infty \\
\left|\arg{\lambda}-\frac{2 \pi k}{5}\right| < \frac{\pi}{5} -\varepsilon}}\frac{\psi_l(\lambda)}{\psi_{l+2}(\lambda)} \in \mathbb{C}
\cup \infty \, , \; k \in \mathbb{Z}_5 \; .
\end{equation}
the k-th asymptotic value.

If $\psi_{l}$ and $\psi_{l+2}$ are linearly independent then $w_k(l,l+2)=w_m(l,l+2)$ if and only if
$\psi_k$ and $\psi_m$ are linearly dependent.

\begin{Def}
Let $E$ be the (open) subset of the $(a,b)$ plane such that $\psi_0(\lambda;a,b)$ and $\psi_{\pm2}(\lambda;a,b)$ are linearly independent
(its complement in the $(a,b)$ plane is the union of two smooth surfaces \cite{eremenko}). 
On $E$ we define the following functions

\begin{eqnarray} \label{def:u2}
u_2(a,b) &=& \frac{w_{2}(0,-2)}{w_{-1}(0,-2)} \\ \label{def:u-2} 
u_2(a,b) &=& \frac{w_{-2}(0,2)}{w_{1}(0,2)} \\
 \label{def:U}
 U(a,b) &=& \left(
\begin{matrix}
u_2(a,b)-1 \\ 
u_{-2}(a,b) -1
\end{matrix} \right) \; .
\end{eqnarray}
All the functions are well defined and holomorphic. Indeed, due to WKB theory we have that $w_{l+1}(l,l+2)$ is always different from $0$ and $\infty$.
\end{Def}

We can characterize the poles of the int\'egrale tritronqu\'e as the zeroes of $U$.
\begin{Thm}[\cite{piwkb}]\label{lem:Uzero}
The point $a \in \bb{C}$ is a pole of the int\'egrale tritronqu\'ee if and only if there exists $b \in \bb{C}$ such that
$(a,b)$ belongs to the domain of $U$ and
$U(a,b)=0$. In other words $\psi_{-1}(\lambda;a,b)$ and $\psi_2{(\lambda;a,b)}$ are linearly dependent and $\psi_{1}(\lambda;a,b)$ and $\psi_{-2}{(\lambda;a,b)}$ are linearly dependent.
\end{Thm}
We remember that the complex number $b$ in previous lemma is the coefficient of the quartic term in the Laurent
expansion of the tritronqu\'ee solution around $a$ (see Section 2.2 in \cite{piwkb}).

\subsection{WKB Theory}

Let $V(\lambda;a,b)$ be the potential of equation (\ref{eqn:schroedinger}).
We call turning point any zero of $V$.
A Stokes line is any curve in the complex $\lambda$ plane along which
the real part of the action is constant, such that at least one turning point belong to its boundary.
The union of all the Stokes line and all turning points is called the Stokes complex of the potential.

\begin{figure}[htbp]
\begin{center}
\input{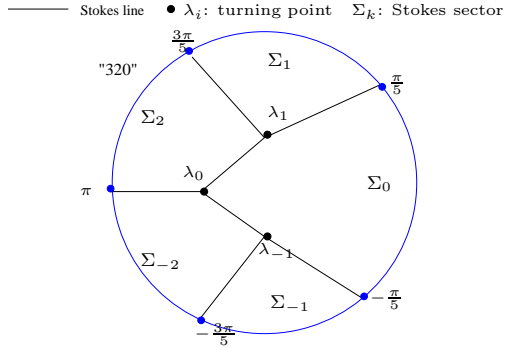}

\end{center}
\caption{Graph "320": dots on the circle represents asymptotic directions in the complex plane}
\label{figure:320}
\end{figure}

A Stokes complex is naturally a graph embedded in the complex plane.
The Stokes graphs has been classified topologically in \cite{piwkb} and the graph of type
"320" (see Figure \ref{figure:320}) was shown to be crucial to
the approximate description of the poles of the int\'egrale tritronqu\'ee.

\begin{Def}
Let $(a^*,b^*)$ be a point such that the Stokes graph of $V(^.;a,b)$ is of type "320". On
a sufficiently small neighborhood of $(a^*,b^*)$
we define the following analytic functions
\begin{eqnarray}\label{def:chi}
 \chi_{\pm2}(a,b) &=& \oint_{c_{\mp 1}}\sqrt{V(\lambda;a,b)}d\lambda \; , \\ \label{def:tu}
 \tilde{u}_{\pm2}(a,b) &=& -e^{\chi_{\pm2}(a,b)} \; ,\\ \label{def:tU}
 \tilde{U}(a,b) &=& \left(
\begin{matrix}
\tilde{u}_2(a,b)-1 \\ 
\tilde{u}_{-2}(a,b) -1
\end{matrix} \right) \; .
\end{eqnarray}

The cycles $c_{\pm1}$ are depicted in Figure 1 and the branch of $\sqrt{V}$ is chosen such that
$\mbox{Re}\sqrt{V(\lambda)} \to + \infty$ as $\lambda \to \infty$ along the positive semi-axis
in the cut plane.

\end{Def}

\begin{Def}\label{def:bsb}
We say that $(a,b)$ satisfies the Bohr-Sommerfeld-Boutroux (B-S-B) system if
the Stokes graph of $V(^.;a,b)$ is of type "320" and
\begin{eqnarray}\nonumber
\chi_{2}(a,b) &=& \oint_{c_{- 1}}\sqrt{V(\lambda;a,b)}d\lambda = i\pi (2n-1) \; , \\ \nonumber
& &~~~~~~~~~~~~~~~~~~~~~~~~~~~~~~~~~~~~~~~~~~~~~~~~(\ref{sys:bsb}) \\ \nonumber
\chi_{-2}(a,b) &=& \oint_{c_{1}}\sqrt{V(\lambda;a,b)}d\lambda = i\pi (2m-1) \; .
\end{eqnarray}
Here $m,n \in \bb{N}-0$ are called quantum numbers.
Equivalently the B-S-B system can be written as $\tilde{U}(a,b)=0$.
\end{Def}

In \cite{kapaev03} the following lemma was proven.

\begin{Lem}\label{lem:uniq}
For any pair of quantum numbers $n,m \in \bb{N}-0$ there exists
one and only one solution of the B-S-B system.
\end{Lem}

After the results of \cite{piwkb} Section 4.3, we can compare the functions $U$ and $\tilde{U}$ defined above.

\begin{Lem}\label{lem:uapprox}
Let $(a,b)$ be such that the Stokes graph is of type "320". There exists a neighborhood of $(a,b)$ and two
continuous positive functions $\rho_{\pm2}$ such that $\chi_{\pm2}$ are holomorphic and
\begin{equation}\label{app:basic}
\av{\tilde{u}_{\pm2} - u_{\pm2}} \leq \frac{1}{2}(e^{2 \rho_{\pm2}}-1) \; .
\end{equation}
Moreover if $\rho_{\pm2} < \frac{\ln 3}{2}$ then $\psi_0$ and $\psi_{\pm2}$ are linearly independent.
\end{Lem}
We remark that in \cite{piwkb} $\rho_{\pm2}$ were denoted $\rho^0_{\pm2}$.

Using classical relations of the theory of elliptic functions we have the following
\begin{Lem}\label{lem:invertibility}
The map $\tilde{U}$ defined in (\ref{def:tU}) is always locally invertible (hence its zeroes are always simple) and
\begin{eqnarray*}
\frac{\partial \chi_{2} }{\partial a}(a,b) \frac{\partial \chi_{-2} }{\partial b}(a,b) -
\frac{\partial \chi_{-2} }{\partial a}(a,b) \frac{\partial \chi_{2} }{\partial b}(a,b)= -28 \pi i  \;.
\end{eqnarray*}

\begin{proof}
On the compactified elliptic curve $\mu^2=V(\lambda;a,b)$, consider the differentials
$\omega_a=-\frac{\lambda d\lambda}{\mu}$ and $\omega_b=-\frac{d\lambda}{\mu}$. 

It is easily seen that $$\frac{\partial \chi_{\pm2} }{\partial a}(a,b)=\oint_{c_{\mp1}}\omega_a \, , \;
\frac{\partial \chi_{\pm2} }{\partial b}(a,b)=14\oint_{c_{\mp1}}\omega_b \; .$$
Moreover we have that
$$J\tilde{U}=\left( \frac{\partial \chi_{2} }{\partial a}(a,b)
\frac{\partial \chi_{-2} }{\partial b}(a,b) - \frac{\partial \chi_{-2} }{\partial a}(a,b)
\frac{\partial \chi_{2} }{\partial b}(a,b)\right)\tilde{u}_{2}\tilde{u}_{-2} \; ,$$
where $J\tilde{U}$ is the Jacobian of the map $\tilde{U}$.

The statement of the lemma follows from the classical Legendre relation between complete elliptic periods
of the first and second kind \cite{bateman2}.
\end{proof}

\end{Lem}

Our aim is to locate the zeroes of $U$ (the poles of the int\'egrale
tritronqu\'ee after Theorem \ref{lem:Uzero})
knowing the location of zeroes of $\tilde{U}$(the solutions of the B-S-B system).
We want to find a neighborhood of a given solution of the B-S-B system
inside which there is one and only one zero of $U$.
Due to estimate (\ref{app:basic}) and Rouch\'e theorem, it is sufficient to find a domain on whose boundary
the following inequality holds
\begin{equation}\label{inequality:rouche}
\frac{1}{2} \left( e^{2\rho_{2}}-1\right)\av{u_2} + \frac{1}{2} \left( e^{2\rho_{-2}}-1\right) \av{u_{-2}} <
 \av{1-\tilde{u}_2}+ \av{1-\tilde{u}_{-2}} \; .
\end{equation}

\subsubsection{Scaling Law}

In order to analyze the important inequality (\ref{inequality:rouche}), we take advantage of the following scaling behaviour that was proven in \cite{piwkb} Section 4.4.

\begin{Lem}\label{lem:chiderivatives}
Let $(a^*,b^*)$ be such that the Stokes graph is of type "320" and $E$ be a neighborhood of $(a^*,b^*)$
such that the estimates (\ref{app:basic}) are satisfied. Then, for any real positive $x$
the point $(x^2a^*, x^3b^*)$ is such that the Stokes graph is of type "320" and in the neighborhood
$E_X=\left\lbrace (x^2a,x^3b): (a,b) \in U \right\rbrace $ the estimates (\ref{app:basic}) are satisfied.
Moreover for any $(a,b) \in E$ the following scaling laws are valid
 \begin{itemize}
  \item $\chi_{\pm2}(x^2a,x^3b)= x^{\frac{5}{2}}\chi_{\pm2}(a,b)$.
  \item $\frac{\partial^{(n+m)} \chi_{\pm2} }{\partial a^n \partial b^m}(x^2a,x^3b)=x^{\frac{5-4 n -6 m}{5}}
  \frac{\partial^{(n+m)} \chi_{\pm2} }{\partial a^n \partial b^m}(a,b)$.
  \item $\rho_{\pm2}(x^2a,x^3b)= x^{-\frac{5}{2}}\rho_{\pm2}(a,b)$.
 \end{itemize}
\end{Lem}

\section{Proof of the main theorem}

From Lemma \ref{lem:chiderivatives} we can extract the leading behaviour of $\tilde{U}$ around real solutions of the
B-S-B system.

\begin{Lem}\label{lem:linearization}
Let  $z=(2k+1)^{\alpha}(a-a^q_k)$, $c_{\pm2}=\frac{\partial \chi_{\pm2} }{\partial a}(a^q,b^q)$, $w=(2k+1)^{\beta}(b-b^q_k)$,  and $d_{\pm2}=\frac{\partial \chi_{\pm2} }{\partial b}(a^q,b^q)$.
If $\alpha > \frac{1}{5}$ and $\beta > -\frac{1}{5}$, then
\begin{eqnarray}\nonumber
\tilde{u}_{2}(z,w)\! &=& \!1+ c_2 (2k+1)^{\frac{1}{5} -\alpha}z + d_2 (2k+1)^{-\frac{1}{5} -\beta}w + O((2k+1)^{-\gamma'}) \; ,
 \\ \label{exp:us} \\ \nonumber 
\tilde{u}_{-2}(z,w) \! &=&\! 1+ c_{-2} (2k+1)^{\frac{1}{5} -\alpha}z +
d_{-2} (2k+1)^{-\frac{1}{5} -\beta}w +
O((2k+1)^{-\gamma'}) \; , \\ \nonumber
&&\gamma' > -\frac{1}{5}+\alpha, \gamma'>\frac{1}{5} +\beta \; . 
\end{eqnarray}

\begin{proof}
It follows from Lemma \ref{lem:chiderivatives}.
\end{proof}
\end{Lem}

\begin{Def}\label{def:polydisc}
We denote $D_{\epsilon, \delta}^{(a',b')}=\left\lbrace\av{a-a'} <\e , \av{b-b'}<\delta, \,
\e,\delta \neq 0\right\rbrace$ the polydisc centered at $(a',b')$.
\end{Def}

We have collected all the elements for proving the following

\begin{Lem}\label{lem:polydisc}
Let $\e$, $\delta$ be arbitrary positive numbers. If $\frac{1}{5} <\alpha < \frac{6}{5}$,
$-\frac{1}{5} <\beta < \frac{4}{5}$, then there exists a $K \in \bb{N}^*$ such that for any $k\geq K$,
$U$ and $\tilde{U}$ are well-defined and holomorphic on $D_{k^{-\alpha}\e,k^{-\beta} \delta}^{(a^q_k,b^q_k)}$
and the following inequality holds true

\begin{equation}\label{ineq:UU}
\av{U(a,b) - \tilde{U}(a,b)} < \av{\tilde{U}(a,b)}, \forall (a,b) \in \partial D_{k^{-\alpha}\e,k^{-\beta} \delta}^{(a_k,b_k)} \; . 
\end{equation}

\begin{proof}
The polydisc $D_{k^{-\alpha}\e,k^{-\beta} \delta}^{(a^q_k,b^q_k)}$ is the
image under rescaling  $a \to (2k+1)^{\frac{4}{5}}a$,
$b \to (2 k+1)^{\frac{6}{5}}b $ of a shrinking polydisc centered at $(a^q,b^q)$; call it $\tilde{D}_k$.
Hence due to Lemma \ref{lem:uapprox}, for $k\geq K'$ $\tilde{D}_k$ is such that
$\rho_{\pm2}$ are bounded, $\chi{\pm2}$ are holomorphic and the estimates
(\ref{app:basic}) hold. Call $\rho^*$ the supremum of
$\rho_{\pm2}$ on  $D_{K'}$. Due to scaling property, for all $k \geq K'$ $\rho_{\pm2}$ is bounded from above by $(2k+1)^{-1}\rho^*$ on $D_{k^{-\alpha}\e,k^{-\beta} \delta}^{(a_k,b_k)}$; such a bound is eventually smaller than $\frac{\ln 3}{2}$.

Then for a sufficiently large $k$, $D_{k^{-\alpha}\e,k^{-\beta} \delta}^{(a_k,b_k)}$ is a subset of the domain of $U$ and inside it $U$ and $\tilde{U}$ satisfy (\ref{app:basic}) and (\ref{exp:us}).

We divide the boundary in two subsets
$\partial D_{k^{-\alpha}\e,k^{-\beta}\delta}^{(a^q_k,b^q_k)}=D_0 \cup D_1$, $$D_0=\left\lbrace |a-a^q_k|=k^{-\alpha}\e; \av{b-b^q_k} \leq k^{-\beta}\delta\right\rbrace \,, $$ $$ D_1 = \left\lbrace |a-a_k| \leq k^{-\alpha}\e; \av{b-b_k} = k^{-\beta}\delta\right\rbrace \; .$$
Inequality (\ref{ineq:UU}) will be analyzed separately on $D_0$ and $D_1$.

If  $\av{d_2}\leq \av{d_{-2}}$, denote $d_2=d, d_{-2}=D, c=c_2, C=c_{-2}$; in the opposite case
$\av{d_2} > \av{d_{-2}}$, denote $d_{-2}=d, d_{2}=D, c=c_{-2}, C=c_{2}$.
By the triangle inequality and expansion (\ref{exp:us}), we have that
$$
\av{\tilde{U}(a,b)} \geq  (2k +1)^{\frac{1}{5} - \alpha} \e \av{(c - \frac{C d}{D})} + \mbox{ higher order terms } , \;
(a,b) \in D_0 .
$$

Similarly, if $\av{c_2}\leq \av{c_{-2}}$ denote $d_2=d, d_{-2}=D, c=c_2, C=c_{-2}$; in the opposite case
$\av{c_2} > \av{c_{-2}}$, denote
$d_{-2}=d, d_{2}=D, c=c_{-2}, C=c_{2}$. By the triangle inequality and expansion (\ref{exp:us}), we have that

$$
\av{\tilde{U}(a,b)} \geq (2k +1)^{-\frac{1}{5} - \beta}\delta \av{(d - \frac{D c}{C})} + \mbox{ higher order terms }, \;
(a,b) \in D_1 .
$$

We observe that $(c - \frac{C d}{D})\neq0$ and $ (d - \frac{D c}{C})\neq0$, since (see Lemma \ref{lem:invertibility}) $c_2d_{-2}- c_{-2}d_2 = - 28 \pi i $. By hypothesis
$ -1<\frac{1}{5} - \alpha<0$ and $ -1<-\frac{1}{5} - \beta<0$.

Conversely, $\av{U(a,b) - \tilde{U}(a,b)} \leq \frac{\rho^*}{2k+1} + \mbox{ higher order terms},$ for all $(a,b) \in D_0\cup D_1$.

The Lemma is proven.

\end{proof}

\end{Lem}

As a corollary of Lemma \ref{lem:polydisc} and of Rouch\'e theorem, we obtain the following theorem which
implies Theorem \ref{thm:aapprox}.

\begin{Thm}\label{thm:approximation}
Let $\e$, $\delta$ be arbitrary positive numbers. If $\frac{1}{5} <\alpha < \frac{6}{5}$,
$-\frac{1}{5} <\beta < \frac{4}{5}$, then it exists a $K \in \bb{N}^*$ such that for any $k\geq K$,
inside the polydisc $D_{k^{-\alpha}\e,k^{-\beta} \delta}^{(a^q_k,b^q_k)}$ there is one and only one solution of
the system $U(a,b)=0$.
\end{Thm}

\paragraph{Acknowledgments}
I am indebted to my advisor Prof. B. Dubrovin who constantly gave me
suggestions and advice.


\bibliographystyle{alpha}

\end{document}